\documentclass[12pt]{article}
\usepackage[dvips]{epsfig}
\usepackage{t1enc,amsmath,amsfonts,amssymb,amscd,amsthm,latexsym}
\usepackage{a4wide}
\usepackage{dcolumn}
\usepackage{nicefrac}
\newcolumntype{d}[1]{D{?}{}{#1}}

\def\tg{\hbox {\rm tan\,}}
\def\ctg{\hbox {\rm ctn\,}}
\def\arctg{\hbox {\rm arctan\,}}

\renewcommand{\Re}{{\rm I\kern-0.16em R}}

\def\@begintheorem#1#2{\trivlist \item[\hskip \labelsep{\bf #1\ #2}]}
\def\@opargbegintheorem#1#2#3{\trivlist
      \item[\hskip \labelsep{\bf #1\ #2\ (#3)}]}

\newtheorem{proposition}{Proposition}[section] 
\newtheorem{defn}[proposition]{Definition}

\newtheorem{thm}[proposition]{Theorem}

\newtheorem{example}[proposition]{Example}
\newtheorem{remark}[proposition]{Remark}

\def\P{{\bf P}}
\def\R{{\bf R}}
\def\R{{\bf R}}
\def\Q{{\bf Q}}
\def\E{{\bf E}}

\def\cE{{\cal E}}
\def\cF{{\cal F}}

\def\cB{{\cal B}}

\def\cG{{\cal G}}

\def\al{\alpha}

\def\Ga{\Gamma}
\def\ga{\gamma}
\def\si{\sigma}

\numberwithin{equation}{section}

\begin{document}

\title{\vspace*{-1.5cm}
\begin{flushright}
\begin{minipage}{4cm}
\tiny
{\sc Dept. of Math. \hfill Univ. of Oslo\\
Pure Mathematics \hfill     No.~35\\
ISSN 0806--2439 \hfill     November 2005}
\end{minipage}
\end{flushright}
\vskip1cm
Perpetual integral functionals of diffusions and
their numerical computations}

\date{\today}

\author{Paavo Salminen\\{\small Åbo Akademi,}
\\{\small Mathematical Department,}
\\{\small Vänriksgatan 3 B,}
\\{\small FIN-20500 Åbo, Finland,}
\\{\small email: phsalmin@abo.fi}
\and
Olli Wallin\\{\small University of Oslo,}
\\{\small Centre for Mathematics and Applications,}
\\{\small P.O. Box 1053 Blindern,}
\\{\small NO-0316 Oslo, Norway}\\
{\small email: olli.wallin@cma.uio.no}}
\vskip5cm

\date{}

\maketitle

\begin{abstract}
In this paper we study perpetual integral functionals of
diffusions. Our interest is focused on cases where such
functionals can be expressed as first hitting times for some other
diffusions. In particular, we generalize the result in
\cite{salminenyor04} in which one-sided functionals of Brownian motion
with drift are connected with first hitting times of reflecting
diffusions.

Interpretating perpetual integral functionals as hitting times
allows us to compute numerically their distributions by applying
numerical algorithms for hitting times. Hereby, we discuss two approaches:
\begin{itemize}
\item numerical inversion of the Laplace transform of the first hitting time,
\item numerical solution of the PDE associated with the distribution
function of the first hitting time.
\end{itemize}
For numerical inversion of Laplace tranforms we have implemented the
Euler algorithm developed by Abate and Whitt. However, perpetuities
lead often to diffusions for which the explicit forms of the Laplace
transforms of first hitting times are not available. In such cases, and
also otherwise,  algorithms for numerical solutions of PDE's can be
evoked. In particular, we analyze the Kolmogorov PDE of some
diffusions appearing in our work via the Crank-Nicolson scheme. 
\\ \\
{\rm Keywords:}  Bessel process, geometric Brownian motion,
random time change, local time.
\\ \\ 
{\rm AMS Classification:} 60J65, 60J60, 62E25.
\end{abstract}

\section{Introduction}
\label{sec0}

Let $\{Y_t\,:\,t\geq 0\}$ be a regular linear diffusion taking values
on an interval $I$. The left and right endpoints of the interval are
denoted by $l$ and $r,$ respectively. For a locally integrable
function $f:I\mapsto \R_+$ define the perpetual integral functional
associated with $f$ and $Y$ via
\begin{equation}
\label{01}
\int_0^\infty f(Y_t)\, dt.
\end{equation}
An important example of perpetual integral functionals is
\begin{equation*}
\label{02}
\int_0^\infty \exp\left(-2a\,B^{(\mu)}_t\right)\, dt, \quad a>0,
\end{equation*}
where $B^{(\mu)}$ is a BM with positive drift $\mu,$
studied by Dufresne in \cite{dufresne90} in connection with risk
theory and pension funding. In particular, from \cite{dufresne90},
this functional is distributed as $1/(2\, a^2\,Z_{\nu})$ where
$Z_{\nu}$ is a gamma-distributed random variable
with the density function
$$
f_{Z_\nu}(z)= \frac{1}{\Gamma(\nu)}\, z^{\nu-1}{\rm e}^{-z}, \quad \nu:=\mu/a.
$$
In Yor \cite{yor92b} (see \cite{yor01} for
an English translation) it is shown
\begin{equation}
\label{d-y}
\int_0^\infty  \exp(-2aB_s^{(\mu)})\,
ds
\quad{\mathop=^{\rm{(d)}}}\quad
H_0(R^{(\delta)}),
\end{equation}
where $R^{(\delta)}$ is a Bessel process of dimension
$\delta=2(1-(\mu/a))$ started at $1/a,$
$$
H_0(R^{(\delta)}):=\inf\{ t:\ R^{(\delta)}_t=0\},
$$
and $\displaystyle{{\mathop=^{\rm{(d)}}}}$ reads
"is identical in law with" (in fact, $R^{(\delta)}$ can be constructed
in the same probability space as $B^{(\mu)}$ and then (\ref{d-y}) holds a.s.).
In \cite{salminenyor04} the methodology used in \cite{yor92b} is
developed for more general perpetual functionals for BM with positive
drift and, in particular, results for one sided functionals are
presented. An example of these is
\begin{equation*}
\label{s-y}
\int_0^\infty  \exp(-2aB_s^{(\mu)})\,{\bf 1}_{\{B^{(\mu)}_s>0\}} ds\,
\quad{\mathop=^{\rm{(d)}}}\quad
H_{1/a}(R^{(2\mu/a)}),
\end{equation*}
where the Bessel diffusion $R^{\,(2\mu/a)}$ is started at 0 and, in
the case $0<\mu<a,$ reflected at 0. For further results and references
for Dufresne's functionals, see \cite{deschepperetal92},
\cite{milevsky97}, \cite{salminenyor03b}, \cite{salminenyor04} and
\cite{decampsetal05}.

In this paper, Section 2, we recall (from \cite{borodinsalminen04}) the connection between
perpetual integral functionals and first hitting times. After this,
the result in \cite{salminenyor04},
Proposition 2.3, concerning one-sided perpetual functionals of  $B^{(\mu)}$
is generalized for $Y$ (defined via a SDE) and functionals of the type in (\ref{01}).
In Section 3, to make the paper more self contained and also as an
introduction to Section 4,
some basic facts about the distributions of the first
hitting times are presented. Section 4 contains brief descriptions of
the Euler algorithm for numerical inversion of Laplace transforms and
the Crank-Nicolson scheme for solving PDE's, which we implemented in Matlab. The paper is concluded
with Section 5  where the distributions of some perpetual
functionals are computed numerically.  In particular, we compare
the one-sided functionals
$$
\int_0^{\infty} \exp(-2B^{(\mu,\sigma)}_s))\,{\bf
1}_{\{B^{(\mu,\sigma)}_s>0\}}\, ds\quad{\rm and}\quad
\int_0^{\infty} (1+\exp(B^{(\mu,\sigma)}_s))^{-2}\,{\bf 1}_{\{B^{(\mu,\sigma)}_s>0\}}\, ds,
$$
where $B^{(\mu,\sigma)}_t:=\sigma\, B_t+\mu\, t$ with
$B$ the standard Brownian motion. It is also seen that some of the diffusions studied have
bad singularities making the PDE's numerically troublesome to solve. In some cases
this problem can, at least partly, be solved by transforming the diffusion to a new one
with better behaviour. It seems to us that for a general numerical
approach 
for calculating distributions of 
perpetualities, more sophisticated PDE or 
other methods such as Monte Carlo simulation are needed for the cases where the Laplace transform 
is not available for numerical inversion.

\section{Perpetual integral functionals as first hitting times}
\label{sec1}

Consider a diffusion $Y$ on an open interval $I=(l,r)$
determined by the SDE
\begin{equation}
\label{a1}
dY_t=\sigma(Y_t)\,dB_t+b(Y_t)\,dt,
\end{equation}
where $B$ is a standard Brownian motion defined in a complete
probability space $(\Omega,\cF,\{\cF_t\},\P).$
It is assumed that $\sigma$ and $b$ are continuous and
$\si(x)>0$ for all $x\in I$. The diffusion $Y$ is considered up to
its explosion (or life) time
$$
\zeta:=\inf\{t\ :\ Y_t\not\in I\}.
$$

Let $f$ be a positive and continuous function defined on $I,$ and
consider for $t\geq 0$ the  integral functional 
\begin{equation*}
\label{A1}
A_t:=\int_0^t f(Y_s)\,ds.
\end{equation*}
We remark that $\{A_t\,:\, t\geq 0\}$ is an additive functional of $Y$
in the usual sense (see e.g. \cite{blumenthalgetoor68}
p. 148). Taking $t=\zeta$  gives us the perpetual integral functional 
\begin{equation*}
\label{A11}
A_\zeta:=\int_0^\zeta f(Y_s)\,ds.
\end{equation*}
Assuming that $A_\zeta<\infty$ a.s. we are interested in the distribution of
 $A_\zeta.$ 

A sufficient condition for finiteness is clearly that 
the mean of $A_\zeta$ is finite:
\begin{eqnarray*}
&&\E_x\left(A_\zeta\right)= 
\int_0^\infty \E_x\left(f(Y_s)\right)\,ds \\
&&\hskip1.5cm
= \int_l^r G_0(x,y)\, f(y)\,m(dy)<\infty,
\end{eqnarray*}
where $G_0$ denotes the Green kernel of $Y$ and $m$ is the speed measure
(for these see, e.g., \cite{borodinsalminen02}).
A neccessary and sufficient condition in the case of a 
Brownian motion with drift $\mu>0$ is that the function $f$ is
integrable at $+\infty$ (see Engelbert and Senf \cite{engelbertsenf91} and 
Salminen and Yor \cite{salminenyor03}). We refer also to a recent
paper \cite{salminenyor06} for such a condition valid for continuous
$f$ and a fairly general diffusion $Y.$  

Next proposition connects the perpetual integral
functionals to the first hitting times. The result is extracted from
Propositions 2.1 and 2.3 in
\cite{borodinsalminen04} where the proof can be found. We remark
also that the result generalizes
Proposition 2.1 in \cite{salminenyor04}.

\begin{proposition}
\label{!}
 {\sl
Let $Y,$ $A,$ and $f$ be as above and assume that there exists a two
times 
continuously differentiable function $g$ such that 
\begin{equation}
\label{fg}
f(x)=\big(g'(x)\si(x)\big)^2,\quad x\in I.
\end{equation}
Let $\{a_t\,:\,0\leq t<A_\zeta\}$ denote the inverse of $A,$ that is,
$$
a_t:=\min\big\{s: A_s>t\big\}, \qquad t\in [0,A_\zeta).
$$
{\bf 1.} Then the process $Z$ given by
\begin{equation}
\label{a2}
Z_t:=g\left(Y_{a_t}\right), \qquad t\in[0,A_\zeta),
\end{equation}
is a diffusion satisfying the SDE
\begin{equation*}
\label{a3} 
dZ_t=d\widetilde B_t+G(g^{-1}(Z_t))\,dt, \qquad t\in[0,A_\zeta).
\end{equation*}
where $\widetilde B_t$ is a Brownian motion and
\begin{equation}
\label{a4} 
G(x)=\frac{1}{f(x)}\left( \frac 12\, \si(x)^2\, g''(x)
+b(x)\,g'(x)\right).
\end{equation}

\noindent
{\bf 2.} Let $x\in I$ and $y\in I$
be such that $\P_x$-a.s.
$$
H_y(Y):=\inf\{t:\, Y_t=y\}<\infty.
$$
Then
\begin{equation*}
\label{aa10}
A_{H_y(Y)}=\inf\{t:\, Z_t=g(y)\}=: H_{g(y)}(Z)\quad a.s.
\end{equation*}
with $Y_0=x$ and $Z_0=g(x).$

\noindent
{\bf 3.} Suppose $g(r):=\lim_{z\to r}g(z)$ exists. 
Suppose also that the following statements hold a.s.
\begin{equation*}
\label{aa100}
(i)\ \ \lim_{t\to\zeta}Y_t=r,\quad 
(ii)\ \ A_\zeta:=\lim_{t\to \zeta}A_t<\infty .
\end{equation*}
Then 
\begin{equation*}
\label{a10}
A_\zeta= H_{g(r)}(Z)\quad a.s.
\end{equation*}
}
\end{proposition}

In \cite{salminenyor04} Proposition 2.3 one sided functionals for Brownian motion with
positive drift are studied. This result is generalized here, under some assumptions,
 to the present case. Suppose $0\in(l,r)$ and recall that
 $f(x)>0$ for all $x\in(l,r).$ Consider the functional 
\begin{equation*}
\label{A2}
A^0_\zeta:=\int_0^\zeta f(Y_s)\,{\bf 1}_{\{Y_s>0\}}\,ds.
\end{equation*}
Let 
$$
C_t:=\int_0^t{\bf 1}_{\{Y_s>0\}}\,ds, \quad t\leq \zeta,
$$
and $\{c_t\,:\,0\leq t<B_\zeta\}$ denote the inverse of $C.$ We assume
also that 
\begin{equation}
\label{A22}
\lim_{t\to\zeta}Y_t=r\quad {\rm  a.s.}
\end{equation}
It is  well known (see \cite{itomckean74}) that the process 
$$
Y^+:=\{Y_{c_t}\,:\,0\leq t<C_\zeta\}
$$  
is identical in law with $Y$ living
on $[0,r)$ and having 0 as a reflecting boundary point. Applying the
  random time change means that on every sample path the excursions below 0 are omitted
  after which the gaps created are closed by joining the excursions
  together. Therefore, 
\begin{equation*}
\label{A3}
A^0_\zeta=\int_0^{\zeta^+} f(Y^+_s)\,ds=:A^+_\zeta
\end{equation*}
where $\zeta^+$ is the life  time of $Y^+.$ 

Next introduce the local time of $Y^+$ at 0  via
\begin{equation*}
\label{L1}
L_t(Y^+):= \sigma^2(0)\,\lim_{\varepsilon\downarrow 0}(2
\varepsilon)^{-1}{\rm Leb}\{0\leq s\leq t\,:\, Y^+_s<\varepsilon\}.
\end{equation*}
Under some additional smoothness assumptions on $\sigma$ and
$b$ (see McKean \cite{mckean69}) the pair $(Y^+,L(Y^+))$ with
$Y^+_0=x>0$ can 
be viewed as the unique solution of the reflected SDE
\begin{equation*}
\label{a11}
dX_t=\sigma(X_t)\,dB_t+b(X_t)\,dt  + dL_t(X),\quad X_0=x,
\end{equation*}
such that 
\begin{description}
\item{(a)} $\lim_{t\to\zeta(X)}X(t)=r,$
\item{(b)} $0\leq X(t)<r$ for all $t<\zeta(X),$
\item{(c)} $t\mapsto L_t(X)$ is continuous, increasing with $L_0(X)=0,$ and 
$$
\int_0^t {\bf 1}_{\{0\}}(X_s)\,dL_s(X)=L_t(X).
$$ 
\end{description}
We now give the promised generalization.
\begin{proposition}
\label{!!}
 {\sl Let $Y^+$ be as given above and define for $t<\zeta^+$ 
$$
A^+_t:=\int_0^{t} f(Y^+_s)\,ds.
$$
The inverse of $A^+$ is denoted by $\{a^+_t\,:\,0\leq t<A_\zeta\}.$
Recall the definition of the function $g$ in (\ref{fg}) and define 
the process $Z^+$ via 
\begin{equation}
\label{a22}
Z^+_t:=g\left(Y^+_{a^+_t}\right), \qquad t\in[0,A^+_\zeta).
\end{equation}
Then 
\begin{equation}
\label{a222}
A^+_\zeta =\inf\{t\,:\, Z^+_t=g(r)\}\quad {\rm a.s.}
\end{equation}
with $Z_0=g(x).$ Moreover, $Z^+$ satisfies the reflected SDE 
\begin{equation}
\label{a33} 
dZ^+_t=d\widetilde B_t+G(g^{-1}(Z^+_t))\,dt + dL_t(Z^+), \qquad t\in[0,A^+_\zeta).
\end{equation}
where $\widetilde B_t$ is a Brownian motion, 
\begin{equation}
\label{a333} 
L_t(Z^+)=\lim_{\varepsilon\downarrow 0}(2
\varepsilon)^{-1}{\rm Leb}\{0\leq s\leq t\,:\, g(0)\leq
Z^+_s<g(0)+\varepsilon\},
\end{equation}
and $G$ is as in (\ref{a4}). The local time $L(Z^+)$ is related to the
local time $L(Y^+)$ by 
\begin{equation}
\label{a3333} 
L_t(Z^+)= g'(0) L_{a^+_t}(Y^+).
\end{equation}
 }
\end{proposition}
\begin{proof}
To fix ideas, assume that $g$ is monotonically increasing. 
By Ito's formula for $u< \zeta$
\begin{eqnarray*}
\label{a7}
\nonumber
&&\hskip-1.5cm
g(Y^+_u)-g(Y^+_0)=\int_0^u g'(Y^+_s)\left(\si(Y^+_s)\,dB_s+b(Y^+_s)\,ds  + dL_s(Y^+)\right)
\\
\nonumber
&&\hskip4cm
+\frac 12 \int_0^u\,g''(Y^+_s)\si^2(Y^+_s)\,ds.
\end{eqnarray*}
Replacing $u$ by $a^+_t$ yields
\begin{eqnarray*}
\label{a71}
\nonumber
&&\hskip-1.5cm
Z^+_t-Z^+_0=\int_0^{a^+_t}g'(Y^+_s)\si(Y^+_s)\,dB_s + g'(0)\,L_{a^+_t}(Y^+)
\\
\nonumber
&&\hskip3cm
+\int_0^{a^+_t}\big(g'(Y^+_s)\si(Y^+_s)\big)^2G(Y^+_s)\,ds.
\end{eqnarray*}
Since $a^+_t$ is the inverse of $A^+_t$ and $(A^+_s)'=\big(g'(Y^+_s)\sigma(Y^+_s)\big)^2$
we have
\begin{equation}
\label{a5}
(a^+_t)'=\frac{1}{(A^+_{a^+_t})'}
=\left(g'(Y^+_{a^+_t})\sigma(Y^+_{a^+_t})\right)^{-2}.
\end{equation}
From L\'evy's theorem it follows that
$$
\widetilde B_t:=\int_0^{a^+_t}g'(Y^+_s)\sigma(Y^+_s)\,dB_s,
\qquad t\in[0,A^+_\zeta),
$$
is a (stopped) Brownian motion. Consequently, for $t<A^+_\zeta$
\begin{eqnarray*}
&&\hskip-.5cm
Z^+_t-Z^+_0=\widetilde B_t+  g'(0)\,L_{a^+_t}(Y^+)+\int_0^t
\big(g'(Y^+_{a^+_s})\si(Y^+_{a^+_s})\big)^2G(Y^+_{a^+_s})\,da^+_s\\
&&\hskip1.3cm
=\widetilde B_t+\int_0^t G(g^{-1}(Z^+_s))\,ds+  g'(0)\,L_{a^+_t}(Y^+).
\end{eqnarray*}
Clearly, viewing  $t\mapsto g'(0)\,L_{a^+_t}(Y^+)$ as a functional of
$Z^+$ then this functional increases only on the set $\{t\,.:\, Z^+_t=g(0)\}.$ Moreover,
since $Y^+_t\geq 0$ for $t\geq 0$ we have $Z^+_t\geq g(0)$ for $t\geq
0$ by monotonicity of $g.$ Hence, $(Z^+,L(Z^+))$ can be seen as the unique solution of 
the reflected SDE (\ref{a33}) with $L(Z^+)$ as in (\ref{a333})
satisfying (\ref{a3333}), as claimed. Finally, again by the
monotonicity of $g,$ the identity (\ref{a222}) follows from the
definition (\ref{a22}) of $Z^+$ and the assumption (\ref{A22}).
\end{proof}

\begin{remark}
Notice that the above approach yields a stronger result than in
\cite{salminenyor04}, i.e., the identity (\ref{a222}) holds a.s. 
\end{remark}

\section{Reminder on first hitting times}
\label{sec2}

\subsection{Distribution functions and PDEs}
\label{sec21}

Let $Y$ be a linear diffusion determined via the SDE (\ref{a1}). It is
here assumed that $Y$ hits $r$ a.s. and is killed when this happens.
Therefore, the boundary point $r$ is either exit-not-entrance or
regular with killing, and $l$ is either natural or entrance-not-exit
or regular with reflection. Letting $H_r(Y)$ denote the
hitting time of $r$ we have 
\begin{equation}
\label{b1}
\P_x(H_r(Y)>t)=\int_l^r  p(t;x,y)\, m(dy),
\end{equation}
where $ p$ denotes the symmetric transition density of $ Y$ with
respect to its speed measure $m.$ It is well known (see
\cite{itomckean74} p. 149 and \cite{mckean56}) 
that $(t,x)\mapsto 
p(t;x,y)$ satisfies for all $y\in (l,r)$ the PDE
\begin{eqnarray*}
\label{b2}
&&
\frac{\partial}{\partial t} p(t;x,y)=\frac 12 \sigma^2(x)\,
\frac{\partial^2}{\partial x^2}
p(t;x,y)+b(x)\,\frac{\partial}{\partial x} p(t;x,y)\\
&&
\nonumber
\hskip2.1cm
=:(\cG\, p)(t;x,y)
\end{eqnarray*}
and the condition
$
\lim_{x\to r} p(t;x,y)=0.
$
Moreover, in the case $l$ is regular with reflection or
entrance-not-exit we impose at $l$ the condition
$$
\lim_{x\to l}\frac{\partial}{\partial x}
p(t;x,y)=0,
$$ 
and in the case $l$ is natural the condition
$$
\lim_{x\to l} p(t;x,y)=\lim_{x\to l}\frac{\partial}{\partial x}
p(t;x,y)=0.
$$
See \cite{itomckean74} or \cite{borodinsalminen02} for the boundary
classification of linear diffusions. 

Letting $\{ T_t\}$ denote the semigroup associated with $ Y$ we
may write from (\ref{b1})
$$
\P_x(H_r( Y)>t)= T_t1(x).
$$
Recall from \cite{mckean56} (where the case with natural scale is
treated) 
that $(t,x)\mapsto  T_tg(x)$ with $g$
bounded and continuous satisfies 
$$
\frac{\partial}{\partial t}( T_tg)(x)=(\cG\, T_t)g(x).
$$
Consequently, the
distribution function 
$$ 
(t,x)\mapsto u(t,x):=\P_x(H_z( Y)<t)
$$
is the unique solution of the PDE problem 
\begin{equation}
\label{b3}
\frac{\partial}{\partial t}u(t,x)=(\cG\,u)(t,x)
\end{equation}
with the initial condition 
$
\lim_{t\to 0}u(t,x)=0\quad {\rm for\ all\ } x\in(l,r)
$
and the boundary condition
$
\lim_{x\to r}u(t,x)=1\quad {\rm for\ all\ } t>0.
$
Further, if $l$ is regular with reflection or
entrance-not-exit 
$$
\lim_{x\to l}\frac{\partial}{\partial x}u(t,x)=0,
$$ 
and in the case $l$ is natural 
$$
\lim_{x\to l} u(t,x)=\lim_{x\to l}\frac{\partial}{\partial x}u(t,x)=0.
$$

\begin{remark}
Using the fact (see \cite{mckean56}) that 
$$
(t,x,y)\mapsto \frac{\partial}{\partial t} p(t;x,y)
$$
is continuous and satisfies the same boundary conditions as the
density $p(t;x,y)$ it is easy (at least when $m(l,r)<\infty  $) to deduce
that 
$$
\frac{\partial}{\partial t}\P_x(H_z( Y)<t)=-
\lim_{y\to z}\frac{1}{S'(y)}\, \frac{\partial}{\partial y} p(t;x,y),
$$
where 
$$
S'(y)=\exp\left(-\int^y 2\,\sigma^{-2}(v)\,b(v)\, dv\right)
$$
is the derivative of the scale function $S$ (cf. \cite{itomckean74} p. 154).
\end{remark}

\subsection{Laplace transforms and ODEs}
\label{sec22}
For the approach with the Laplace transform of $H_r(Y)$ consider
 the second order ODE
\begin{equation}
\label{f1}
\cG u(x)=\lambda u(x), 
\end{equation}
where $\lambda\geq 0.$ It is known (see \cite{feller52} p. 488,
and \cite{itomckean74} p. 128) that the equation (\ref{f1}) has a
positive increasing solution $\psi_\lambda$ and a
positive decreasing solution $\varphi_\lambda.$ In case $l$ is natural
or entrance and $r$ is exit these solutions are unique up to 
multiplicative constants. When $l$ is regular with reflection
the condition $\psi'_\lambda(l)=0$ must be posed, and 
when $r$ is
regular with killing the condition is $\varphi_\lambda(r-)=0.$   
The Green kernel $G_\lambda$ of $Y$ can be expressed via these
solutions  as
\begin{eqnarray*}
&&
G_\lambda(x,y):=\int_0^\infty {\rm e}^ {-\lambda\,t}\, p(t;x,y)\, dt\\
&&
\hskip1.8cm
=
\begin{cases}
\frac{1}{w_\lambda}\ \psi_\lambda(x)\,\varphi_\lambda(y),&\ x\leq y,\\
\frac{1}{w_\lambda}\ \psi_\lambda(y)\,\varphi_\lambda(x),&\ y\leq x,\\
\end{cases}
\end{eqnarray*}
where $w_\lambda$ is the Wronskian (see e.g. \cite{borodinsalminen02}).
Using the Green kernel the Laplace transform for the first hitting
time $H_y(Y)$ is given by 
$$
\E_x\left( {\rm e}^{-\lambda\, H_y(Y)}\right)={\displaystyle
  \frac{G_\lambda(x,y)}{G_\lambda(y,y)}}
$$
and, in particular,  
\begin{equation}
\label{f2}
\E_x\left( {\rm e}^{-\lambda\, H_r(Y)}\right)={\displaystyle
  \frac{\psi_\lambda(x)}{\psi_\lambda(r)}}.
\end{equation}

\section{Numerical methods}
\label{sec3}
\subsection{Numerical inversion of Laplace transforms}
\label{ssec31}
There are several efficient methods for numerical inversion of the Laplace
transforms 
of probability density functions or probability distribution functions.
We have implemented a method developed by Abate and Whitt in
\cite{abatewhitt95b}. This so called the Euler-algorithm has proved to
be very effective in many applications see
e.g. \cite{fumadanwang97}, \cite{carrlinetsky00}, \cite{davydovlinetsky01} and 
\cite{davydovlinetsky01b}.
The main features of
this method are presented below. For more details and also for
further references, see    
\cite{abatewhitt95b}.

Consider a non-negative  random variable
with density $f$ and its Laplace transform 
$$
\hat {f}(\lambda):=\int_0^\infty {\rm e}^{-\lambda t}\, f(t)\,dt.
$$
The well known inversion integral formula (called the Bromwich or also
the Fourier-Mellin integral) states that
\begin{equation}
\label{mellin}
f(t) = \frac{1}{2\pi\, {\rm i}}\int_{a-{\rm i}\,\infty}^{a+{\rm i}\,\infty} {\rm
e}^{\lambda\,t}\, \hat {f}(\lambda)\, d\lambda,
\end{equation}
where it is assumed that $\hat f$ does not have singularities on
or to the right of the vertical line $\lambda =a.$ 
\begin{remark}
\label{spectral}
For first hitting times of diffusions considered in Section 3 we have 
(cf. (\ref{f2}))
$$
\hat{f}(\lambda)=\E_x\left( {\rm e}^{-\lambda\, H_r(Y)}\right)=
  \frac{\psi_\lambda(x)}{\psi_\lambda(r)}.
$$
If the left boundary point $l$ is not natural it follows from
the classical theory of second order differential operators 
(see e.g. \cite{kent80}) that ${\psi_\lambda(x)}$ is for every
$x\in(l,r)$ 
an entire function of $\lambda$ and the zeroes of $\lambda\mapsto{\psi_\lambda(x)}$
are  for every
$x\in(l,r)$ simple and negative. Consequently, the inversion formula
(\ref{mellin}) holds in this case for any $a>0.$ If $l$ is natural we
can approximate the first hitting time $H_r(Y)$ via a sequence of
first hitting times $\{H_r(Y^{(n)})\}$ associated with the diffusions
$Y^{(n)},$ $n=1,2,\dots,$ constructed from $Y$ by reflection at
$l+\frac 1n,$ respectively. Then $H_r(Y^{(n)})\to H_r(Y)$ in
distribution and by dominated convergence it is seen that
(\ref{mellin}) is valid also in this case. 
\end{remark}

Since ${\rm Re}(\hat f(a+{\rm i}\,u))={\rm Re}(\hat f(a-{\rm i}\,u)),$ ${\rm Im}(\hat
f(a+{\rm i}\,u))=-{\rm Im}(\hat f(a-{\rm i}\,u)),$ and 
 $f(t)=0$ for $t<0$ the inversion integral (\ref{mellin}) takes the
form 
\begin{equation}
\label{mellin2}
f(t) = \frac{2{\rm e}^{at}}{\pi }\int_{0}^{\infty} 
{\rm Re}(\hat f(a+{\rm i}\,u))\,\cos(ut)\,du.
\end{equation}
Next we approximate $f(t)$ by using 
the trapezoidal rule for the integral 
(\ref{mellin2}) (see \cite{abatewhitt95b} for some comments on
the effectiveness of this procedure). Letting $h$ denote the step size we have for fixed $a$ and $t$ 
$$
f(t) \approx f_h(t) = \frac{h{\rm e}^{at}}{\pi}\,{\rm Re}\big(\hat{f}(a)\big) 
+ \frac{2h{\rm e}^{at}}{\pi}\sum_{k=1}^{\infty}{\rm Re}\big(\hat{f}(a+kh\,{\rm i})\big)\cos(kht).
$$
Choosing $h=\pi/(2t)$, $a=A/(2t)$ (with $A$ to be made precise later) and 
truncating the infinite series to the first $j$ terms we are led to
define 
\begin{equation}
\label{sj}
s_j(t) := \frac{{\rm e}^{A/2}}{t}\sum_{k=0}^j (-1)^k a_k(t),
\end{equation}
where
$$
a_0(t) := \hat f\left(A/2t\right)/2 
$$
and 
$$
a_k(t) := {\rm Re}\left(\hat{f}\left((A+2k\pi\, {\rm
i})/2t\right)\right),
\quad k=1,2,\dots, j.
$$

It is  possible to accelerate the convergence of the series in (\ref{sj})
by considering it as an alternating series (which is not the case in
general) and using the Euler
summation with binomial weights (see \cite{abatewhitt95b}). Hence, the
proposed final approximation with parameters $m,$ $n,$ and $A$ is
\begin{equation*}
\label{euler}
f(t) \approx E(m,n,t) := \sum_{k=0}^m \binom{m}{k}\,2^{-m}\,s_{n+k}(t),
\end{equation*}
In the examples below we use, following \cite{abatewhitt95b}, $m=11$
and $n=15.$ The error associated with Euler summation can be estimated
by considering the difference
$$
  E(m,n+1,t)-E(m,n,t).
$$

It is advantageous from numerical computational point of view to
invert, 
instead of the density function,  the complementary distribution
function. Therefore consider 
$$
\hat {F}^c(\lambda):=\int_0^\infty {\rm e}^{-\lambda t}\, (1-F(t))\,dt,
$$
where $F$ is the distribution function associated
with $f.$ Firstly, the fact $|1-F(t)|\leq
1$ can be used to show (see \cite{abatewhitt95b}) that
$$
|e_d|\leq \frac {{\rm e}^{-A}}{1-{\rm e}^{-A}},
$$
where $e_d$ stands for the discretization error  
when approximating the integral in (\ref{mellin2}) for $\hat F^c$ via
the trapezoidal rule. For instance, $A=18.4$ gives the upper bound $10^{-8}.$ 
Secondly, under some additional smoothness assumption, it can be proved
(see \cite{abatewhitt95b} Remark 1) that for $\hat F^c$ we have
$a_k(t)>0$ when $k/t$ is large enough motivating the use
of the Euler summation (since the
series in (\ref{sj}) is now alternating). 

For the first hitting time of $r$ for the diffusion $Y$   
the Laplace transform of the complementary distribution function is
given by
\begin{eqnarray*}
&&
\hat {F}^c(\lambda)=\frac 1\lambda\left(1-
\E_x\left( {\rm e}^{-\lambda\, H_r(Y)}\right)\right)=\frac 1\lambda\left(1-{\displaystyle
  \frac{\psi_\lambda(x)}{\psi_\lambda(r)}}\right)\\
&&\hskip1.2cm
=\frac 1\lambda\,\frac{\psi_\lambda(r)-\psi_\lambda(x)}{\psi_\lambda(r)}
\end{eqnarray*}

\subsection{Numerical solutions of PDEs}
\label{ssec32}
In this section we describe using \cite{thomas95} and  \cite{strikwerda89} 
two finite difference methods
known as the
\emph{Crank-Nicolson} (C-N) scheme and the \emph{backward Euler} (BE) method. The C-N scheme is used with satisfactory
results in \cite{poulsen} for calculation of transition probability densities of 
certain diffusions. 
The BE method can be applied in connection with the C-N scheme for the first time step to damp some numerical oscillations typical to the C-N scheme.
Both methods are unconditionally stable: the numerical solutions are well behaved (do not blow up) for any choice of $\Delta t$. Because of the 
singularities appearing in the drift coefficients of our equations, methods that do not have this property (such as the explicit, forward Euler method) 
are practically unusable since they would require extremely small time steps.

To start with, let us introduce a uniformly spaced grid on the
rectangle $[0,T]\times [l,r]$ with $(M+1)\times (N+1)$ nodes, that is, for
\[ \Delta t = \frac{T}{M}, \quad \Delta x = \frac{r-l}{N}, \]
we let $t_m=m\Delta t$ and $x_n=l+n\Delta x$ where $m=0,1,...,M$, $n=0,1,...,N$.
For a real valued function $f$ we use the following finite difference approximations, which can
be justified with Taylor's expansion:
\begin{itemize}
\item[$\cdot$]the \emph{forward difference} approximation of the derivative of $f$ is
\[ \frac{\partial f}{\partial x}(x_i) = \frac{f(x_{i+1})-f(x_{i})}{\Delta x} + \mathcal{O}(\Delta x). \]
\item[$\cdot$]
the \emph{centralized difference} approximation of the derivative of $f$ is
\[ \frac{\partial f}{\partial x}(x_i) = \frac{f(x_{i+1})-f(x_{i-1})}{2 \Delta x} + \mathcal{O}((\Delta x)^2). \]
\item[$\cdot$] the centralized difference approximation of the second order derivative of $f$ is
\[ \frac{\partial^2 f}{\partial x^2}(x_i) = \frac{f(x_{i+1})-2f(x_i)+f(x_{i-1})}{(\Delta x)^2} + \mathcal{O}((\Delta x)^2). \]
\end{itemize}
The C-N scheme approximates the left hand side in equation $(\ref{b3})$ with
the forward difference and the right hand side with the average of the centralized differences at two consequtive times. Denoting
$u_n^m := u(t_m,x_n)$ and dropping the truncation error terms, the discretized equation then reads

\begin{eqnarray*}
\frac{u_n^{m+1}-u_n^m}{\Delta t} &=& \frac{1}{2}\sigma^2(x_n)\frac{1}{2}\Big( \frac{u_{n+1}^{m+1}-2u_n^{m+1}+u_{n-1}^{m+1}}{(\Delta x)^2}
+ \frac{u_{n+1}^{m}-2u_n^{m}+u_{n-1}^{m}}{(\Delta x)^2}\Big)\\
&+& b(x_n)\frac{1}{2}\Big( \frac{u_{n+1}^{m+1}-u_{n-1}^{m+1}}{2\Delta x}
+ \frac{u_{n+1}^{m}-u_{n-1}^{m}}{2\Delta x}\Big).
\end{eqnarray*}
Multiply both sides with $\Delta t$, and define $r_1 = \frac{\Delta t}{2\Delta x}$, $r_2 = \frac{\Delta t}{2(\Delta x)^2}$.
Rearranging the terms so that the values at time $t_{m+1}$ appear on the left hand side and values at time $t_{m}$ appear on the right hand side, we have

\[ A_n u_{n-1}^{m+1} + B_n u_{n}^{m+1} - C_n u_{n+1}^{m+1} = -A_n u_{n-1}^m + D_n u_{n}^m + C_n u_{n+1}^m \]
where
\begin{eqnarray}\label{coefa}
A_n &=& \frac{1}{2}\Big(b(x_n)r_1 - \sigma^2(x_n)r_2\Big), \\ \label{coefb}
B_n &=& 1 + \sigma^2(x_n)r_2, \\ \label{coefc}
C_n &=& \frac{1}{2}\Big(b(x_n)r_1 + \sigma^2(x_n)r_2\Big), \\
D_n &=& 1 - \sigma^2(x_n)r_2.
\end{eqnarray}
Together with the boundary conditions, these form a set of $N+1$ linear equations which we then solve for each $m=1,2,...,M+1$, using the initial
condition for $m=0$.
The Neumann boundary condition (which is needed at a reflecting or an entrance boundary point) is implemented with the second order approximation. In other words,
from
\[ \frac{u_1^m-u_{-1}^m}{2\Delta x} = 0\]
we have $u_{-1}^m = u_1^m$ for the value of $u$ at a "ghost" point beyond the boundary. Plugging this into the discretized equation at
$n=0$ gives
\[ (1+r_2\sigma^2(x_n))u_0^{m+1}-r_2\sigma^2(x_n)u_1^{m+1} = (1-r_2\sigma^2(x_n))u_0^m+r_2\sigma^2(x_n)u_1^m. \]
Using a second order approximation for the boundary condition seems to be important in order to have good convergence.\\
For the backward Euler method, one similarly takes the central approximations for the spatial derivatives but now only at the time step $m+1$. This leads to
the equation

\[ A_n u_{n-1}^{m+1} + B_n u_{n}^{m+1} - C_n u_{n+1}^{m+1} = u_{n}^m \]
for $n=1,...,N$ and $m=1,...,M$, where $A_n$, $B_n$, $C_n$ are given in (\ref{coefa})-(\ref{coefc}) but now with $r_1=\frac{\Delta t}{\Delta x}$, 
$r_2 = \frac{\Delta t}{(\Delta x)^2}$. The second order implementation of the Neumann boundary condition becomes
\[ (1+r_2\sigma^2(x_n))u_0^{m+1}-r_2\sigma^2(x_n)u_1^{m+1} = u_{0}^m. \]

Both methods described above are second order accurate in space, but only the C-N scheme is second order accurate in time. 
However, C-N is known to produce numerical oscillations around discontinuities and sharp gradients if the drift term is 
large compared to the diffusion coefficient (see for example \cite{britz}, \cite{gwharrison}), while the BE method does not have this problem. 
In our examples this is seen as oscillations near the killing boundary. Although these oscillations were damped quite rapidly both in space and time, 
they still produced a small phase shift in the numerical solution of the hitting time distribution $t\mapsto u(x,t)$. 
For this reason the first step in the C-N scheme is divided into $10-100$ substeps, as suggested in \cite{britz}. 
For the substeps we used the BE method. While small oscillations still remained for some of the examples, this procedure 
provided sufficient damping to give very accurate results in our test cases.

A problem more serious than the oscillations appearing in the C-N
scheme is
 faced when the diffusing particle is pushed away from an exit boundary 
by a drift term tending to infinity  in the vicinity of the boundary. In such cases  
convergence can be very slow or even nonattainable without huge computer capacity. 
We discuss this problem in more detail in the examples below.

\section{Examples}
\label{sec4}
In this section some perpetuities are examined numerically. If 
the Laplace transforms of the functionals are available we use the Euler
algorithm for computing the density  and/or the distribution functions.
Applying in these cases also the numerical methods based on
the associated 
PDE's obtained from the  hitting 
time representations of the functionals and checking that the solutions
resulting from the different methods coincide we are able to 
verify the correctness of the
implementations. 

Numerical methods for solving PDEs constitute a powerful tool for computing 
hitting time distributions of diffusions in general. However, there
are diffusions as in Example \ref{bessel} for which the methods
presented here work unsatisfactorily. At least in some particular cases 
it is possible to transform the diffusion to a new one for which the
methods seem to work better. This is discussed in Example \ref{bessel}. 
Due to these difficulties one may consider the numerical inversion 
of the Laplace transforms when available as the first choice for
the kind of numerical computations studied in the paper. 

Below the notation $\{ B^{(\mu,\sigma)}_t\,:\,
t\geq 0\}$ is used for a Brownian motion with (infinitesimal) drift
$\mu$ and variance $\sigma,$ i.e., $ B^{(\mu,\sigma)}_t=\sigma\,
B_t+\mu\, t$ with $B$ a standard Brownian motion.  We assume that 
$\sigma>0$ and $\mu>0.$  In
case $\sigma=1$ we write $ B^{(\mu)}$ for $ B^{(\mu,1)}.$ For a Bessel
process 
with dimension parameter
$\delta$ we use the notation of $\{ R^{(\delta)}_t\,:\,
t\geq 0\},$ and refer to \cite{borodinsalminen02} for their properties.
If nothing else is
written it is assumed that $ B^{(\mu,\sigma)}_0=0$ and 
$ R^{(\delta)}_0=0.$

\begin{example}
\label{hypo}
{\rm Consider the perpetual 
integral functional
\begin{equation*}
\label{x1}
I_1:=\int_0^\infty\cosh^{-2}(B^{(\mu)}_t)\, dt.
\end{equation*}
The Laplace transform of this functional is computed in
\cite{borodinsalminen04} and \cite{vagurina04} for an arbitrary
initial value $x;$ taking therein $x=0$ gives
\begin{eqnarray*}
&&
\E_0\left(\exp\Bigl(-\rho\,I_1\Bigr)\right)=K\,_2F_1(\al,\beta,1+\mu;1/2).
\end{eqnarray*}
where 
$$
\alpha=\frac 12 +\frac 12\sqrt{1-8\rho},\quad
\beta=\frac 12 -\frac 12\sqrt{1-8\rho},
$$
$$
K=
\frac{\Gamma\left(\mu+\al\right)\ \Gamma\left(\mu+\beta\right)}
{\Gamma(\mu)\ \Gamma(\mu+1)},
$$
and $_2F_1$ denotes 
Gauss's hypergeometric function given by 
\begin{eqnarray*}
&&\hskip-1.2cm
_2F_1(a,b,c;x)
:=\frac{\Gamma(c)}{\Gamma(a)\,\Gamma(b)}\sum\limits_{k=0}^{\infty}
 \frac{\Gamma(a+k)\,\Gamma(b+k)}{\Gamma(c+k)}\,\frac{x^k}{k!}\\
&&\hskip1.3cm
=1+\sum\limits_{k=1}^{\infty}
 \frac{a(a+1)\dots(a+k-1)\,b(b+1)\dots(b+k-1)}{c(c+1)\dots(c+k-1)}
\frac{x^k}{k!}.
\end{eqnarray*}
 
Applying Proposition \ref{!} with $g(x):=2\arctg {\rm e}^x$ we obtain
\begin{equation*}
\label{d1}
I_1=H_\pi(Z) \qquad a.s.,
\end{equation*}
where $Z$ satisfies 
\begin{equation}
\label{d2}
dZ_t=dB_t+\left(\frac12\,\ctg Z_t+\frac{\mu}{\sin Z_t}\right)\,dt,\qquad
Z_0=\pi/2.
\end{equation}

For the drift term in (\ref{d2}), it holds
$$
G(g^{-1}(x))=\frac12\, \ctg x+\frac{\mu}{\sin x}
=\frac12\,\left(\mu-\frac12\right)\,\tg\frac x2+
\frac12\left(\mu+\frac12\right)\ctg\frac x2.
$$
Notice that for $0<\mu<1/2$ the drift of $Z$ tends to $-\infty$ when
$Z$ is approaching $\pi.$

In Figures \ref{cosh_de} and \ref{cosh_di} we present
the density and the distribution functions, respectively, of
$I_1$ computed with the Euler algorithm. It has been checked that the PDE method
yields the same results. However, the convergence in the case $\mu=0.3$ seems to be slow. In fact, for
$\mu=0.1$ the convergence rate of the PDE method is so slow that we were unable to get
satisfactory results with our limited computation capacity (RAM).\\


\begin{figure}
\label{cosh_de}
\scalebox{0.7}{\includegraphics{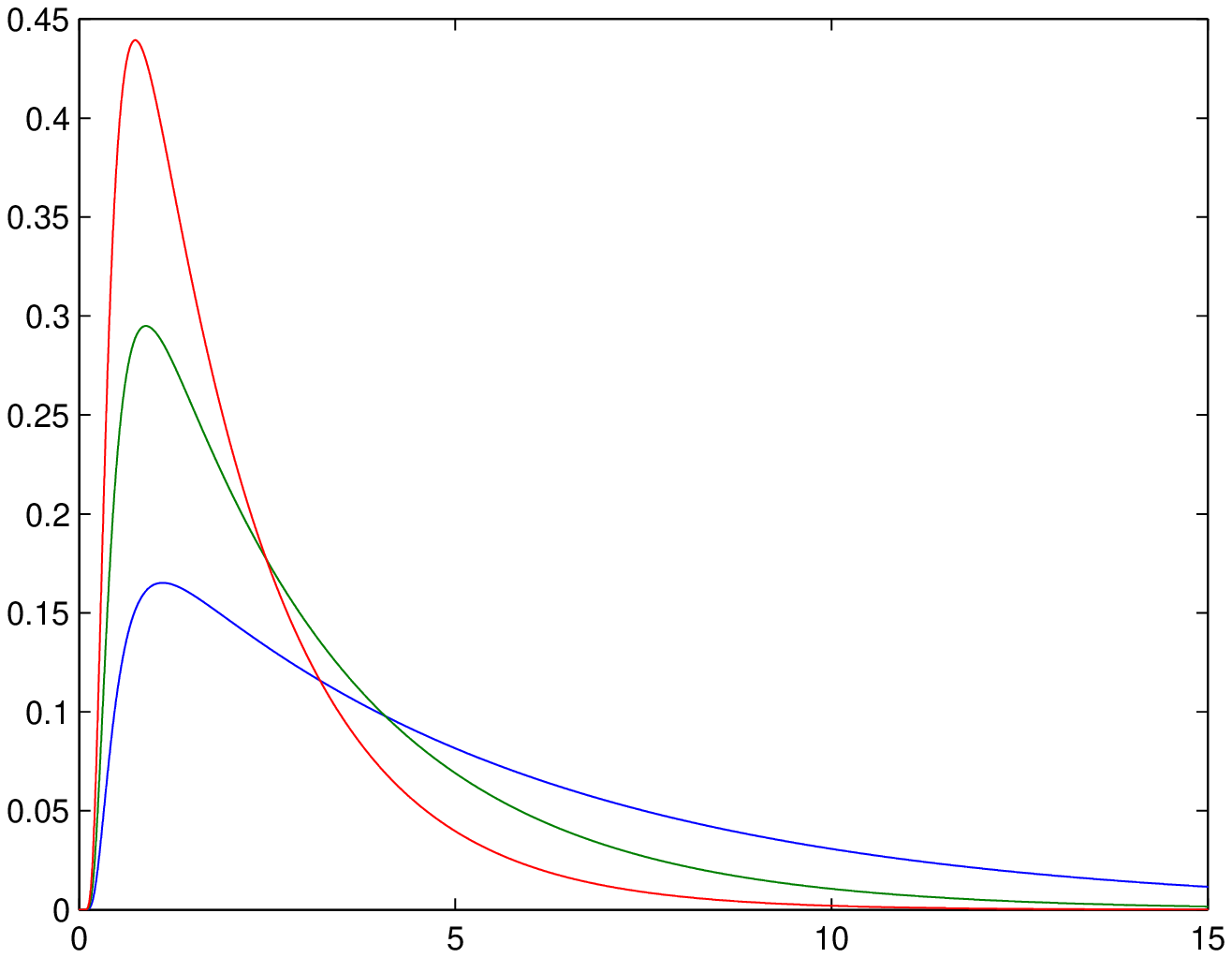}}
\caption{Density of $\int_0^\infty
\cosh^{-2}(B^{(\mu)}_t)\, dt$ for $\mu = 0.3$ (lowest peak),
$\mu=0.5$ and $\mu=0.7$ (highest peak).}

\label{cosh_di}
\scalebox{0.7}{\includegraphics{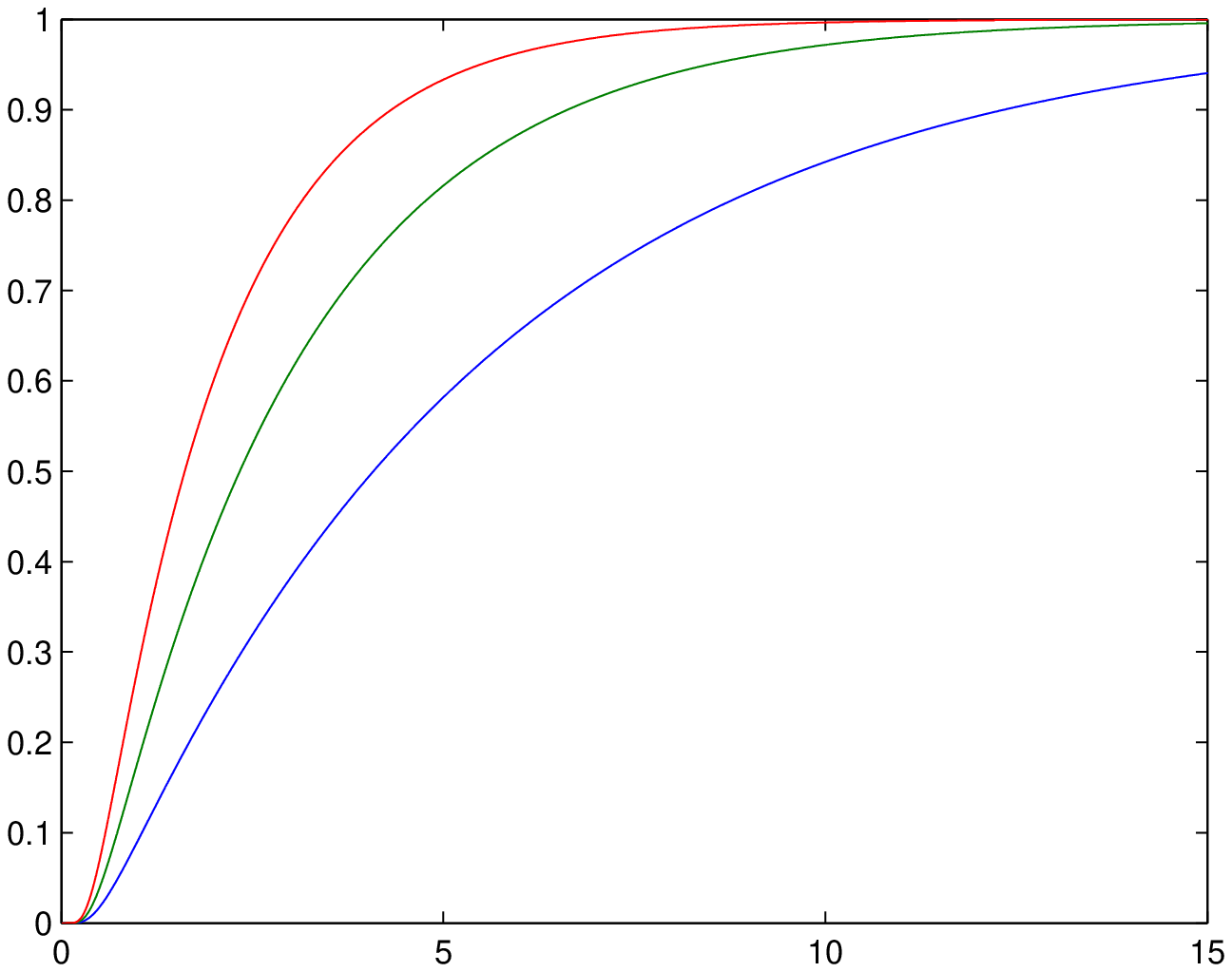}}
\caption{Distribution function of $\int_0^\infty
\cosh^{-2}(B^{(\mu)}_t)\, dt$ for $\mu = 0.3$ (lowest curve),
$\mu=0.5$ and $\mu=0.7$ (highest curve).}
\end{figure}

}
\end{example}

\begin{example}
\label{dufresne}
{\rm In this example we compare the functionals 
$$
I_2:= \int_0^\infty \exp(-2 B^{(\mu,\sigma)}_s) ds\quad{\rm and}\quad 
I_3:= \int_0^\infty ( \exp(B^{(\mu,\sigma)}_s) +1)^{-2}  ds.
$$
Notice that 
$$
I_3= \int_0^\infty  \exp(-2 B^{(\mu,\sigma)}_s)\,\frac 1{ (1+\exp(-B^{(\mu,\sigma)}_s))^2}\, ds,
$$
hence $I_3$ may be seen as a modification of $I_2$ which does
not allow arbitrary large positive discounting. We remark also that
$I_3$ has all moments which is not the case with $I_2.$ 

The Dufresne-Yor identity (cf. (\ref{d-y})) states that 
\begin{equation*}
\label{d3}
I_2 = H_0(R^{(2-2\mu/{\sigma^2})}) \qquad {\rm a.s.},
\end{equation*}
where $R_0^{(2-2\mu/\sigma^2)} = 1/\sigma$. Consequently,  
\begin{eqnarray*}
&&\E_0\left(\exp\bigl(-\rho\,I_2\bigr)\right)
= \frac{\varphi_\rho(1/\sigma)}{\varphi_\rho(0)},
\end{eqnarray*}
with (see \cite{borodinsalminen02} p. 133)
$$
\varphi_\rho(x)= x^{-\nu}\,K_{\nu}(x\sqrt{2\rho}),\quad {\rm and}\quad
\varphi_\rho(0)=2^{-(\nu+2)/2}\Gamma(-\nu)\,\rho^{\nu/2}
$$
and $\nu=-\mu/\sigma^2.$ 

For $I_3$ we have the identity 
\begin{equation*}
\label{d4}
I_3 = H_0(Z)\qquad {\rm a.s.},
\end{equation*}
where $Z$ is the diffusion associated with the SDE
$$
dZ_t=dB_t+\left(\mu +(\mu-\frac{1}{2})\frac{\exp(Z_t)}{1-\exp(Z_t)}\right)\,dt,\qquad
Z_0=-\log 2.
$$
Notice that here $g(x) := -\log(1+\exp(-x))$ (cf. Proposition \ref{!})
and that $Z$ lives on $\R_-$. From  \cite{borodinsalminen04} we recall
the Laplace transform 
\begin{eqnarray*}
&&\E_0\Bigr(\exp\Bigr(-\rho\,I_3
\Bigr)\Bigr)
=
K\
2^{\,\mu-\sqrt{\mu^2+2\rho}}\
_2F_1(\alpha,\beta,\al+\beta +2\mu\,;1/2),
\end{eqnarray*}
where 
$$
\alpha=\frac 12 -\mu +\sqrt{\mu^2+2\rho}+
\sqrt{\frac 14+{2\rho}},\quad
\beta=\frac 12 -\mu +\sqrt{\mu^2+{2\rho}}-
\sqrt{\frac 14+{2\rho}},
$$
and
$$
K=
\frac{\Gamma\left(2\mu+\al\right)\ \Gamma\left(2\mu+\beta\right)}
{\Gamma(2\mu+\al+\beta)\ \Gamma(2\mu)}.
$$
See Figures \ref{duf_de}, \ref{duf_di}, \ref{modduf_de} and \ref{modduf_di}
for illustrations of the distributions of $I_2$ and $I_3$ computed with
the Euler algorithm.

For both functionals in this example it was possible to solve the corresponding
PDE numerically for $\mu\geq\frac{1}{2}$. For $\mu<1/2$ the drift term
tends to $-\infty$ as $Z$ approaches the killing boundary 0. 
This again leads to very slow convergence. While it was 
still possible to achieve good results for some choices of $\mu<\frac{1}{2}$, for a small enough $\mu$ the results were bad even with 
the finest grid we could run on the computer. Notice that here we also need to truncate the semi-infinite domains into finite ones for 
numerical computations. This did not constitute a major problem, but
with larger domains it is difficult to achieve (depending on the
computer capacity) a grid which is spatially dense enough for accurate computations.

\begin{figure}
\scalebox{.70}{\includegraphics{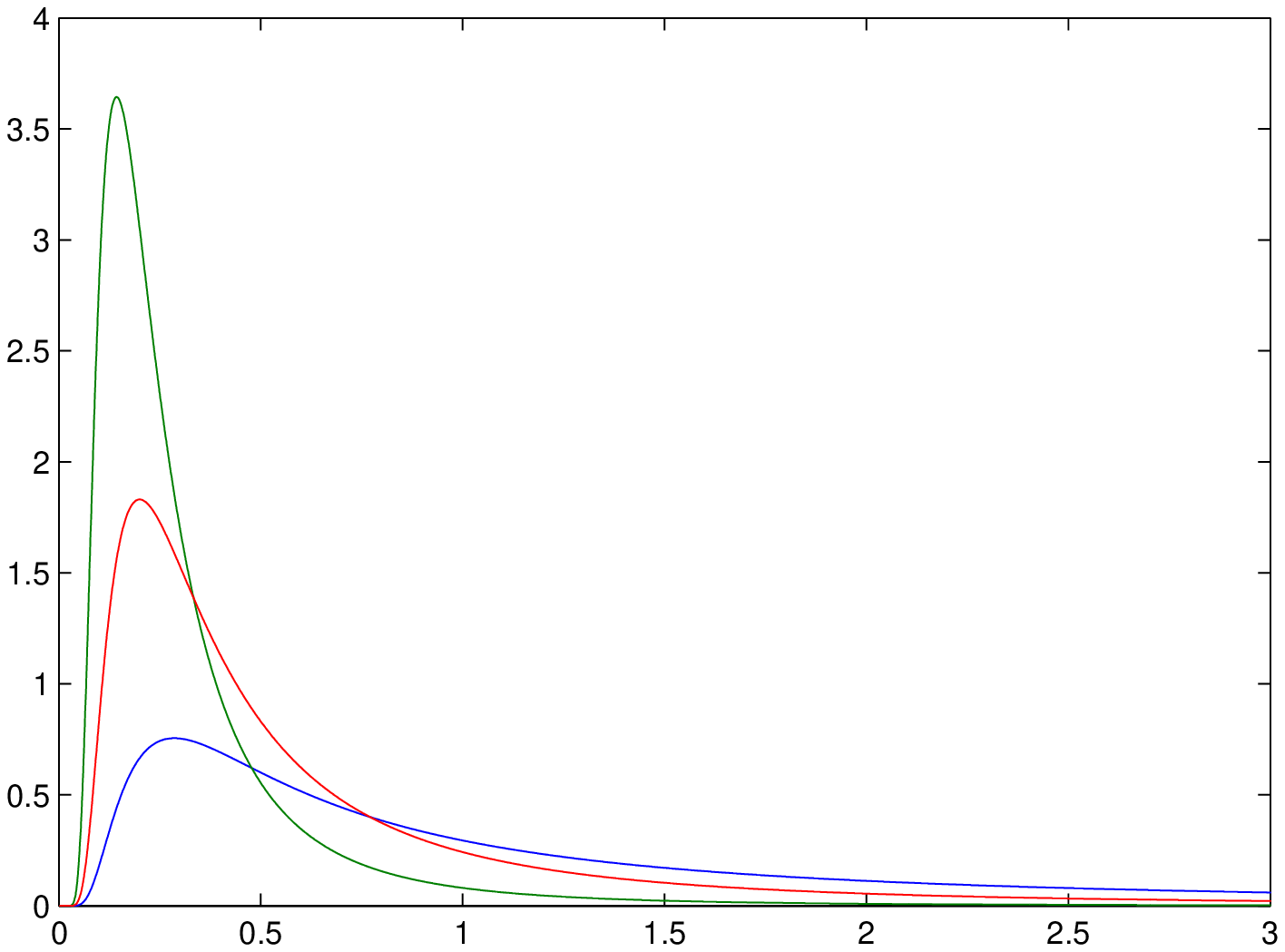}}
\caption{Density of $\int_0^\infty
\exp(-2\,B^{(\mu)}_t)\, dt$ for $\mu = 0.75$ (lowest peak),
$\mu=1.50$ and $\mu=2.50$ (highest peak).}\label{duf_de}

\scalebox{.70}{\includegraphics{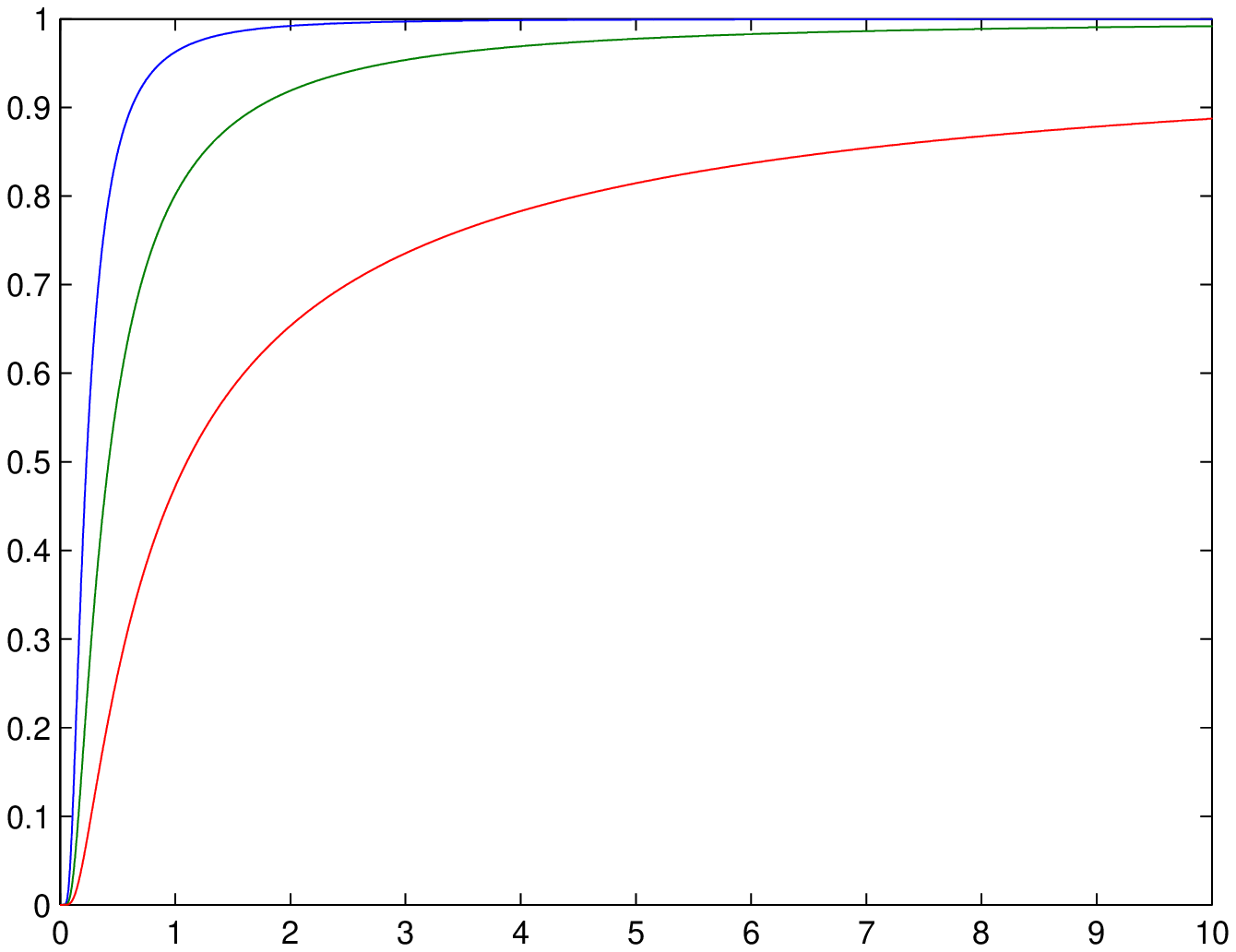}}
\caption{Distribution function of $\int_0^\infty
\exp(-2\,B^{(\mu)}_t)\, dt$ for $\mu = 0.75$ (lowest curve),
$\mu=1.50$ and $\mu=2.50$ (highest curve).}\label{duf_di}
\end{figure}

\begin{figure}
\scalebox{.50}{\includegraphics{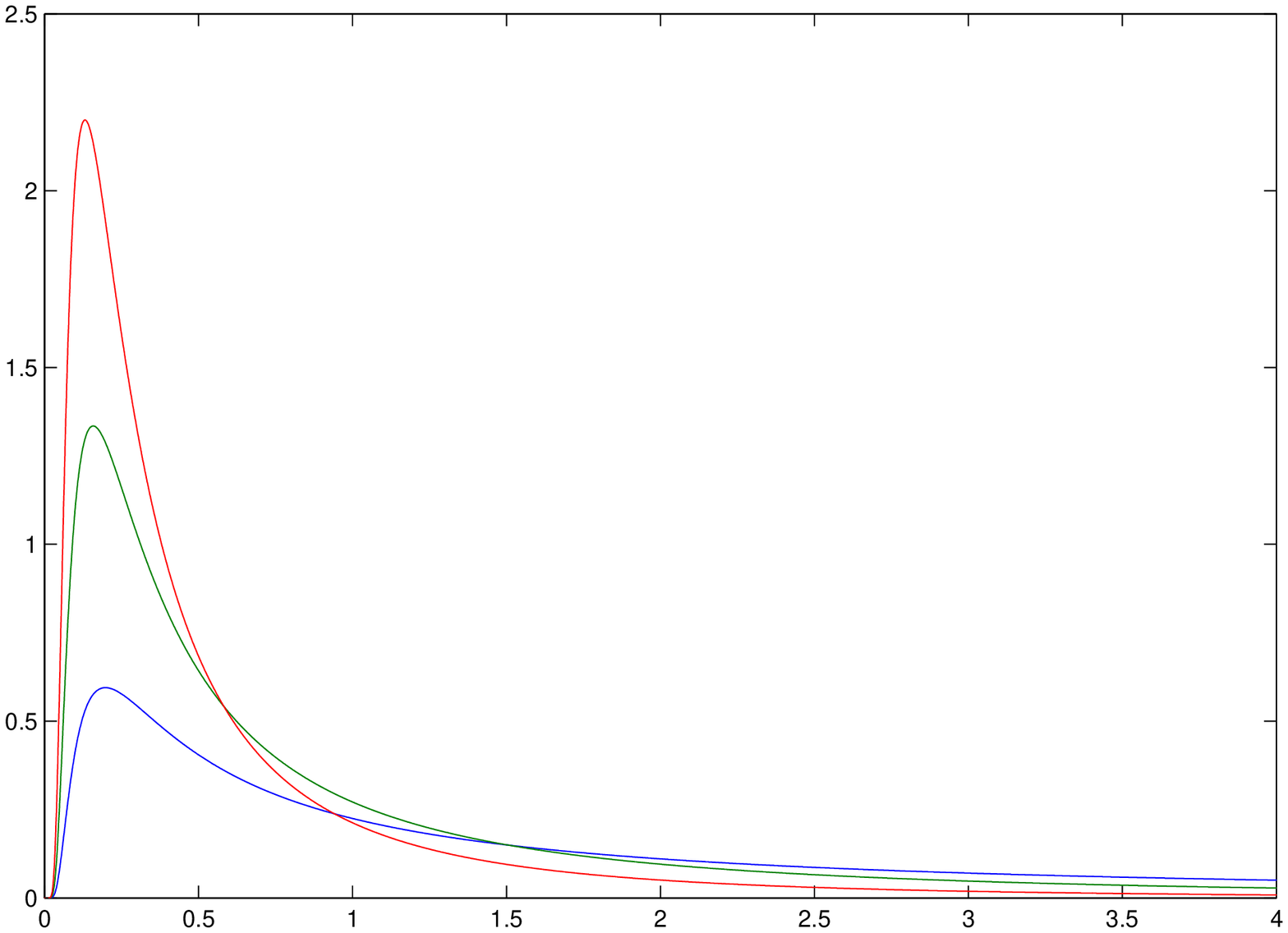}}
\caption{Density function of $\int_0^\infty
(\exp(\,B^{(\mu)}_t)+1)^{-2}\, dt$ for $\mu = 0.25$ (lowest peak),
$\mu=0.50$ and $\mu=0.75$ (highest peak).}\label{modduf_de}

\scalebox{.70}{\includegraphics{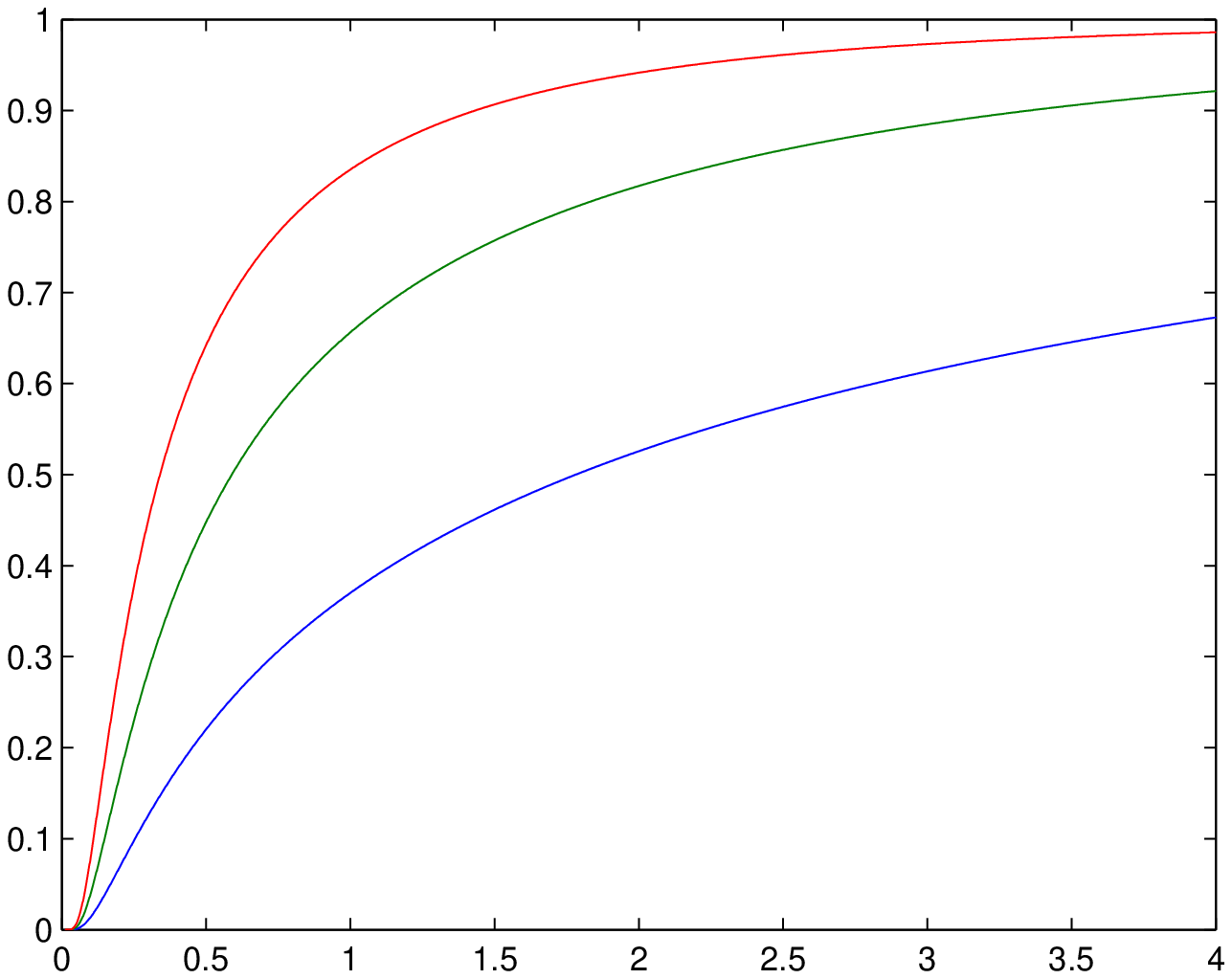}}
\caption{Distribution function of $\int_0^\infty
(\exp(\,B^{(\mu)}_t)+1)^{-2}\, dt$ for $\mu = 0.25$ (lowest curve),
$\mu=0.50$ and $\mu=0.75$ (highest curve).}\label{modduf_di}
\end{figure}

}
\end{example}

\begin{example}
\label{onesided}
{\rm
We define the one-sided variants of $I_2$ and $I_3$ via 
$$
I_4:=\int_0^\infty \exp(-2 B^{(\mu,\sigma)}_s) {\bf
1}_{\{B^{(\mu,\sigma)}_s>0\}} ds 
$$
and 
$$
I_5:=\int_0^\infty ( \exp(B^{(\mu,\sigma)}_s) +1)^{-2} 
{\bf 1}_{\{B^{(\mu,\sigma)}_s>0\}} ds,
$$
respectively.

In \cite{salminenyor04} it is shown that 
\begin{equation*}
\label{d5}
I_4
= H_{1/\sigma}(R^{(2\mu/\sigma^2)})\qquad {\rm a.s.}
\end{equation*}
where $R_0^{(2\mu/\sigma^2)} = 0$. The Laplace transform of $I_4$ is hence given
by
\begin{eqnarray*}
&&\E_0\left(\exp\Bigl(-\rho\,I_4\Bigr)\right)
= \frac{\psi_\rho(0)}{\psi_\rho(1/\sigma)},
\end{eqnarray*}
with (see \cite{borodinsalminen02} p. 133)
$$
\psi_\rho(x)= x^{-\nu}\,I_{\nu}(x\sqrt{2\rho})\quad {\rm and}\quad
\psi_\rho(0)=\frac{\rho^{\nu/2}}{2^{\nu/2}\,\Gamma(\nu+1)}
$$
and $\nu=\mu/\sigma^2-1.$

The Laplace transform of the functional $I_5$  (in \cite{salminenyor04} this
is called the one-sided translated Dufresne functional)
is not known but the following identity (see \cite{salminenyor04}) holds
$$
I_5
= H_0(Z)\qquad {\rm a.s.},
$$
where $Z$ is a diffusion associated with the generator 
$$
\cG f(x)=\frac 12\frac{d^2f}{dx^2}(x)+
\left(\frac 12\sigma +
\frac{\mu-\frac 12\sigma^2}{\sigma\left(1-\exp(\sigma x)\right)}\right)
\frac{d\,f}{dx}(x)
$$
living on $[-(\log 2)/\sigma,0),$  having $-(\log 2)/\sigma$ as a reflecting barrier,
and 0 as a killing barrier.  

\begin{figure}
\scalebox{.60}{\includegraphics{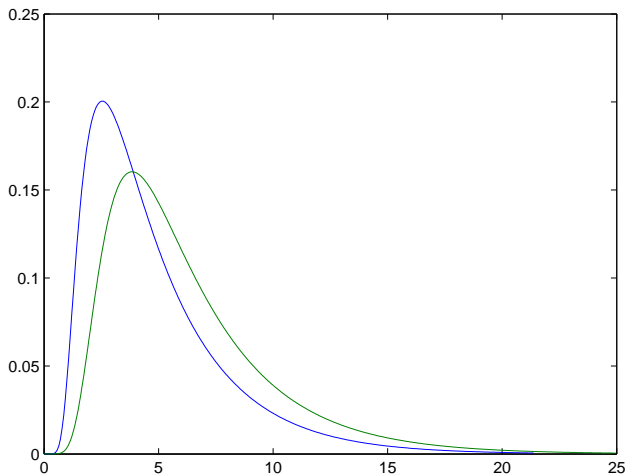}}
\caption{Density function of $\int_0^\infty
(\exp(\,B^{(\mu,\sigma)}_t)+1)^{-2}\, {\bf
    1}_{\{B^{(\mu,\sigma)}_t>0\}}\, dt$ (upper peak) compared with the
  density function of $\int_0^\infty
\exp(-2\,B^{(\mu,\sigma)}_t)\, {\bf
    1}_{\{B^{(\mu,\sigma)}_t>0\}}\, dt$ for  $\mu = 0.04$ and $\sigma=0.20.$ }\label{d2d4}
\end{figure}

In Figure \ref{d2d4} we compare the densities of $I_4$ and $I_5$ 
(see \cite{decampsetal05} for comparisions between $I_2$ and $I_4$). The
density and distribution functions of $I_5$ are displayed for
different values of $\mu$ and $\sigma$ in Figures  \ref{transones_de}
and \ref{transones_di}.

\begin{figure}
\scalebox{.70}{\includegraphics{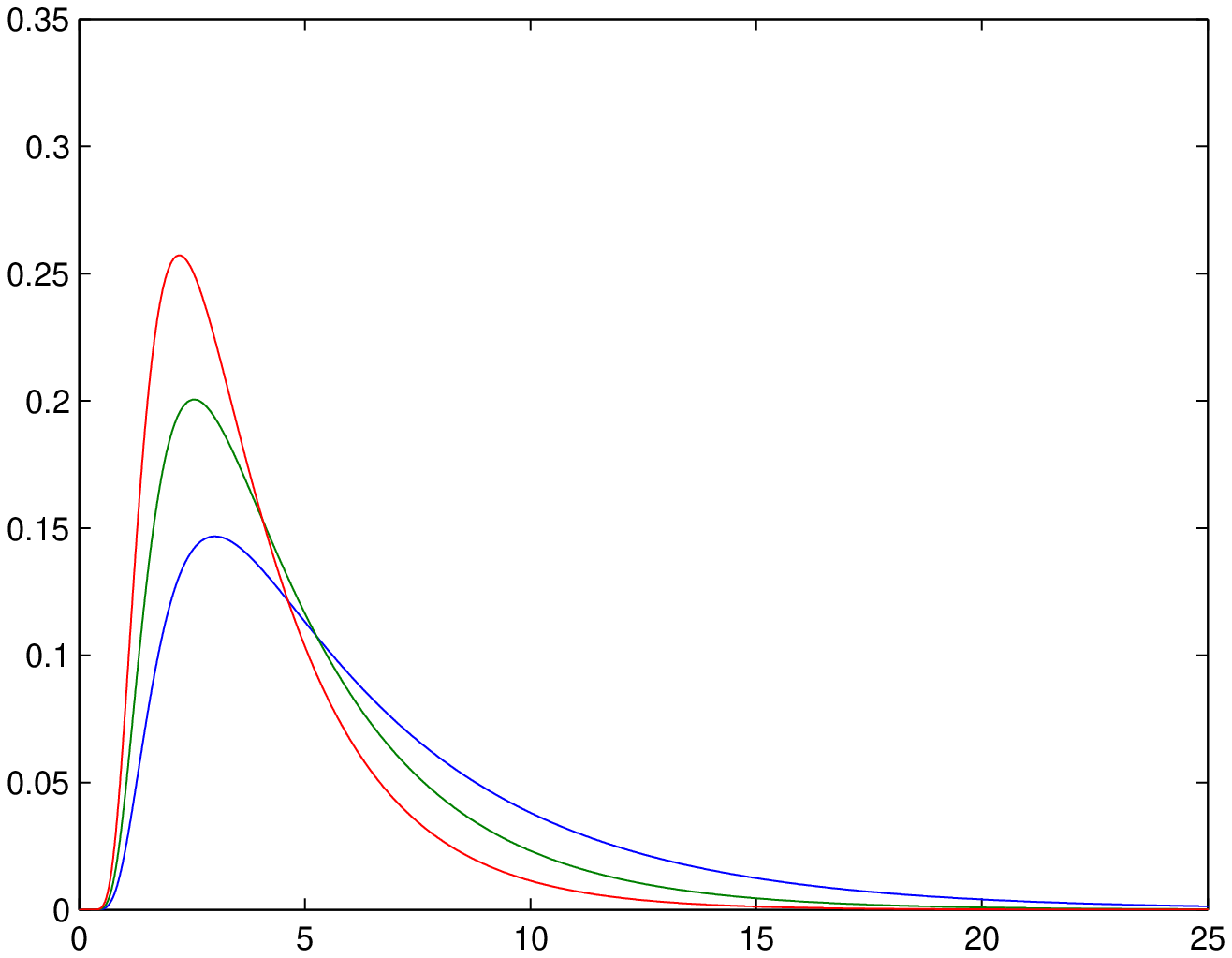}}
\caption{Density function of $\int_0^\infty
(\exp(\,B^{(\mu,\sigma)}_t)+1)^{-2}\, {\bf
    1}_{\{B^{(\mu,\sigma)}_t>0\}}\, dt$ for $\sigma=0.20$ and $\mu =
  0.03$ (lowest peak), $\mu=0.04,$ and
  $\mu=0.05$ (highest peak)}\label{transones_de}

\scalebox{.70}{\includegraphics{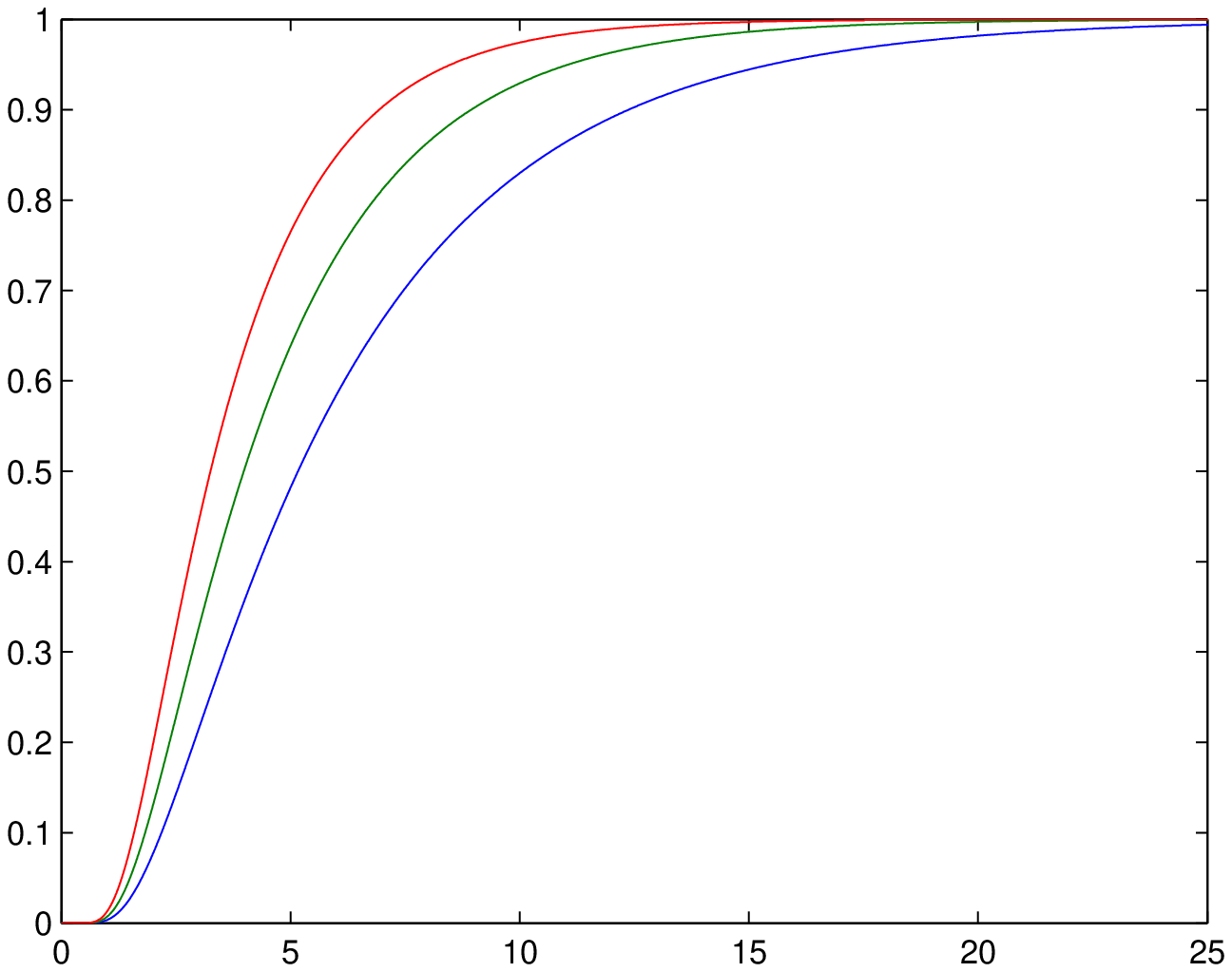}}
\caption{Distribution function of $\int_0^\infty
(\exp(\,B^{(\mu,\sigma)}_t)+1)^{-2}\, {\bf
    1}_{\{B^{(\mu,\sigma)}_t>0\}}\, dt$ for $\sigma=0.20$ and $\mu =
  0.03$ (lowest curve), $\mu=0.04,$ and
  $\mu=0.05$ (highest curve)}\label{transones_di}
\end{figure}

}
\end{example}

\begin{example}
\label{bessel}
{\rm In our final example we consider the functional 
$$
I^{(\delta)}_6:=
\int_0 ^\infty \exp(-2\,R^{(\delta)}_s) ds,\quad  \delta\geq 2.
$$ 
Proposition \ref{!} when applied for $R^{(\delta)}$ and $
g(x) := \exp (x)
$ leads us to the identity
\begin{equation}
\label{d6}
I_6 = H_{0}(Z)\qquad
{\rm a.s.}
\end{equation}
with $Z$ a diffusion associated with 
the SDE 
\begin{equation}
\label{odrift}
dZ_t=dB_t +\frac1{2\,Z_t}\left(1+\frac{\delta-1}{\log Z_t}\right)\,dt,\qquad
Z_0=1.
\end{equation}

In the case $\delta =3$ it is known (see Legall \cite{legall85} and also 
\cite{salminenyor04}) that 
\begin{equation}
\label{d61}
I^{(3)}_6 = H_{1}(R^{(2)})\qquad
{\rm a.s.}
\end{equation}
with $R^{(2)}_0=0.$ 

Since we do not have an
expression for the Laplace transform of $I^{(\delta)}_6$ for
$\delta\not= 3$ we solve numerically the associated
PDE. Unfortunately, 
due to the
complexity of the drift term (in particular, notice that this 
tends, for all values on
$\delta\geq 2$, to $+\infty$  in the vicinity of $0+$) 
simple finite difference schemes do not seem to give solutions converging to the 
correct one, see Figure \ref{d54_huono}.
In search for improvement we implemented a nonuniform grid making the spatial discretization 
denser near the boundaries, and used a fourth-order implementation at
the Neumann boundary. 
While this yielded better results than what is 
seen in Figure \ref{d54_huono}, full convergence still remained out of reach.

\begin{figure}

\scalebox{.60}{\includegraphics{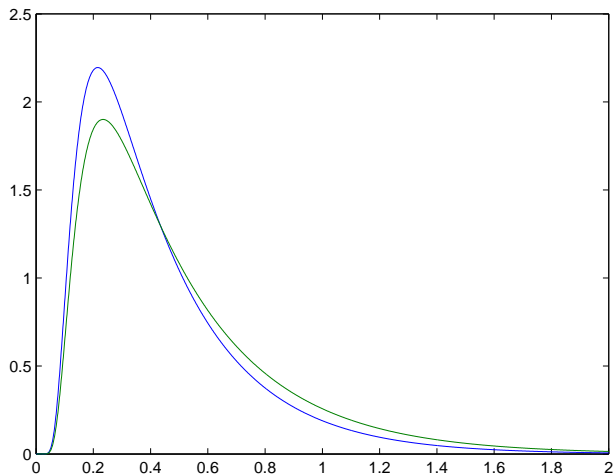}}
\caption{Density function of $\int_0^\infty
\exp(-2\,R^{(\delta)}_t)\, dt$ for $\delta = 3.0$ (upper peak)
obtained using the drift $\frac 1{2x}(1+\frac{\delta-1}{\log x})$
compared with the correct one (lower peak).}\label{d54_huono}
\end{figure}

  These difficulties can at least partly be overcome by 
transforming the diffusion $Z$ given via SDE (\ref{odrift}). 
Indeed, we study now the $h$-transform of $Z$ with $h(x)=S(x) -S(0) $ where $S$
is the scale function of $Z$  Straightforward computations (cf. 
\cite{borodinsalminen02} p. 17) show that we may take 
$$
S(x)=\frac 1{\delta-2} \left|\log x\right|^{2-\delta},\quad 0<x<1
$$
(for simplicity we consider only the case $\delta >2$). 
Then 
$$
\lim_{x\to 0} S(x)=0\quad {\rm and}\quad S'(x)= x^{-1}|\log
x|^{1-\delta}.
$$
Consequently, the generator of the $h$-transform is given by
\begin{eqnarray*}
&&
\cG^\uparrow f=\frac 12 \, \frac{d^2f}{dx^2}+
\frac1{2x}\left(1+\frac{\delta-1}{\log x}\right)\, \frac{df}{dx} + \frac{S'(x)}{S(x)}\, \frac{df}{dx}
\\
&&\hskip.8cm 
=\frac 12 \, \frac{d^2f}{dx^2}+
\frac1{2x}\left(1+\frac{3-\delta}{\log x}\right)\, \frac{df}{dx},\qquad 0<x<1.
\end{eqnarray*}
Let $Z^\uparrow$ denote the $h$-transform, i.e., $Z^\uparrow$ is the
diffusion associated with the generator $\cG^\uparrow.$ By Williams \cite{williams74}
time reversal result (see \cite{borodinsalminen02} p. 35, also for further references)
$$
H_0(Z)\quad{\mathop=^{\rm{(d)}}}\quad H_1(Z^\uparrow).
$$

The PDE associated with $Z^\uparrow$ seems to be well suited for numerical computations. Notice, 
in particular, that if $\delta>3$ the drift term of $Z^\uparrow$ tends to $+\infty$ as $x\to 1$ which fact is in strong contrast  with the corresponding
behaviour of the drift term of $Z$. Hereby it is also of interest to
classify the boundaries of  $Z$ and  $Z^\uparrow.$ It holds for $Z$
that the boundary point $0$ is exit-not-entrance 
and $1$ is entrance-not-exit. For the process $Z^\uparrow$ we have that $0$ is entrance-not-exit and $1$ is entrance-exit (regular) if 
$2<\delta<4$ and entrance-not-exit if $\delta\geq 4$. 
Figures \ref{besduf_de}, \ref{besduf_di} show the density and distribution functions of $I_6$ for some
choises of $\delta$ computed from the PDE associated with $Z^\uparrow.$

As a final comment, and as an extra bonus from our transformation, we
remark that when $\delta=3$ then $\cG^\uparrow$ is the generator of
$R^{(2)},$ and we have recovered the identity (\ref{d61}).

\begin{figure}
\scalebox{.70}{\includegraphics{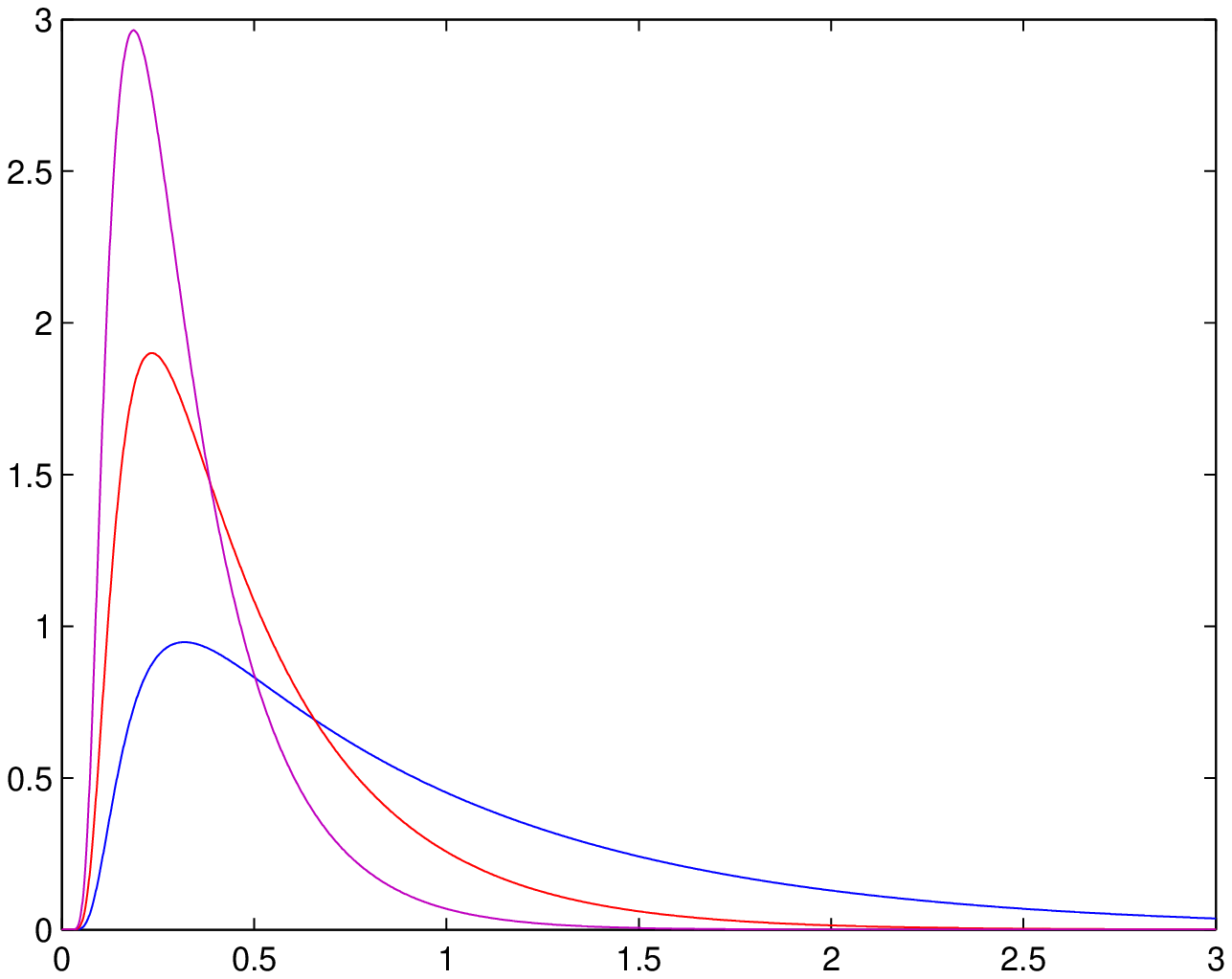}}
\caption{Density function of $\int_0^\infty
\exp(-2\,R^{(\delta)}_t)\, dt$ for $\delta = 2.5$ (lowest peak),
 $\delta = 3.0,$  $\delta = 3.5$ (highest peak).}\label{besduf_de}

\scalebox{.70}{\includegraphics{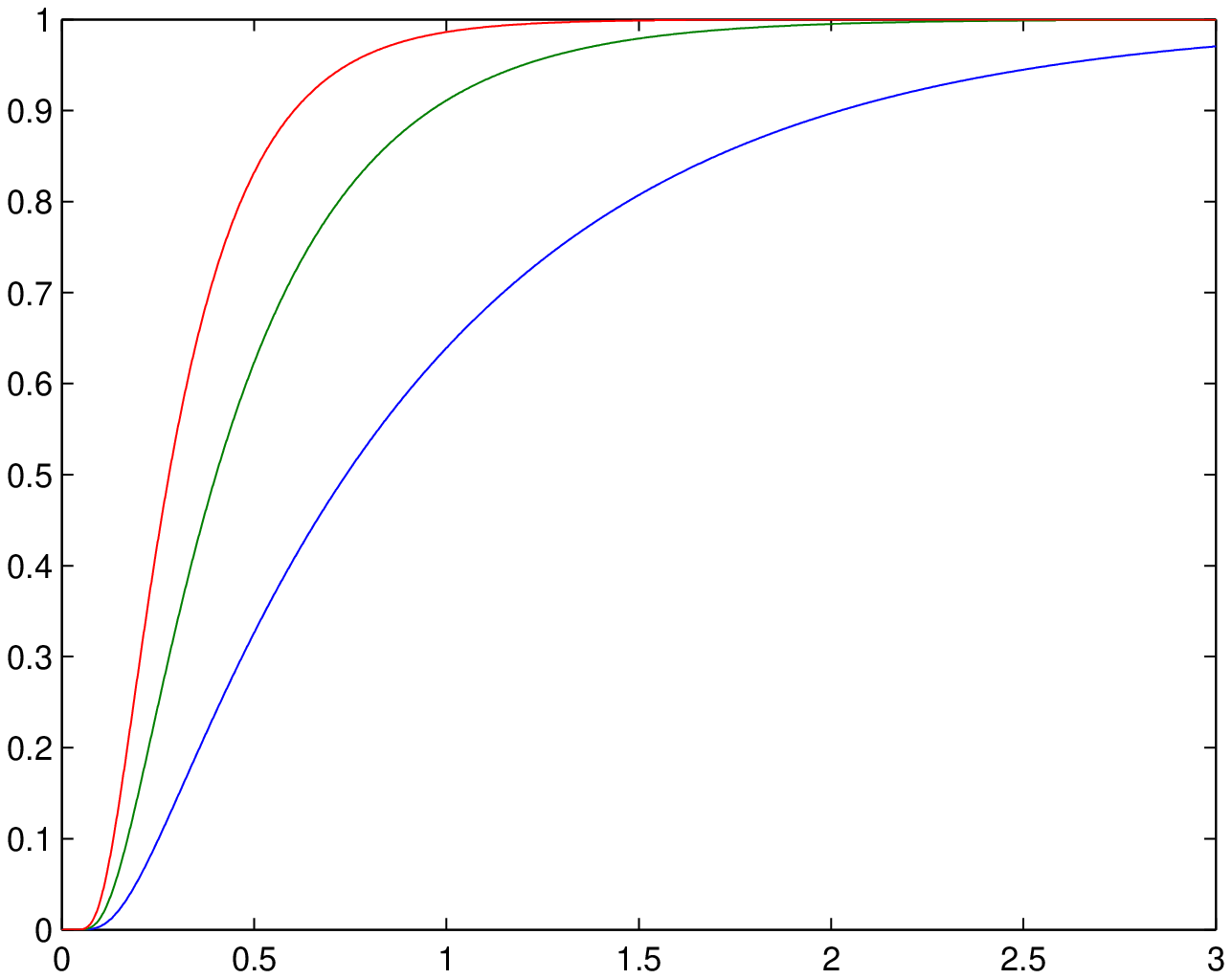}}
\caption{Distribution function of $\int_0^\infty
\exp(-2\,R^{(\delta)}_t)\, dt$ for $\delta = 2.5$ (lowest curve),
 $\delta = 3.0,$  $\delta = 3.5$ (highest curve).}\label{besduf_di}
\end{figure}

}
\end{example}

\textbf{Acknowledgements.} Paavo Salminen thanks Vadim Linetsky for information on numerical
inversion of Laplace transforms of hitting times. Olli Wallin thanks Siddhartha Mishra for several useful 
discussions on numerical solution of partial differential equations. 

\eject
\bibliographystyle{plain}
\bibliography{sal_wal_o}
\end{document}